\documentclass[12pt]{elsart}
\usepackage[latin1]{inputenc}
\usepackage{amsmath}
\usepackage{amsfonts}
\usepackage{amssymb}

\newcommand{\minus}{\!-\!}
\newtheorem{theorem}{Theorem}
\newtheorem{lemma}{Lemma}

\begin{document}
\begin{frontmatter}
\title{Equiangular Tight Frames From Paley Tournaments}
\author{Joseph M.~Renes}
\address{Department of Physics and Astronomy\\
University of New Mexico\\
800 Yale Blvd NE\\
Albuquerque, NM 87131-1156 USA}
\ead{renes@phys.unm.edu}

\begin{abstract}
We prove the existence of equiangular tight frames having $n=2d-1$ elements
drawn from either $\mathbb{C}^d$ or $\mathbb{C}^{d-1}$ whenever $n$ is either $2^k-1$ for $k\in\mathbb{N}$,
or a power of a prime such that $n\equiv 3\!\!\mod 4$. 
We also find a simple explicit expression for the prime power case by establishing
a connection to a $2d$-element equiangular tight frame based on quadratic residues.
\end{abstract}

\end{frontmatter}

\section{Introduction}
\noindent Define an $(n,d)$ equiangular tight frame to be a set of $n$ normalized vectors 
$\phi_k\in\mathbb{C}^d$ such that
\begin{equation}
\label{eq:etf}
|\langle\phi_j,\phi_k\rangle|^2=\frac{n\minus d}{d(n\minus 1)}\qquad \forall\;\, j\neq k\,.
\end{equation}
This is merely the equiangular criterion, but any frame meeting this condition is automatically tight~\cite{bf03,sh03}. 
Many constructions of equiangular tight frames are known~\cite{sh03}, from the ordinary regular simplex, to frames based on 
perfect difference sets~\cite{Koenig03}, quadratic residues~\cite{zauner99}, conference matrices~\cite{ls73,chs96}, 
Hadamard matrices~\cite{hp04}, adjacency matrices of strongly regular graphs or regular two-graphs~\cite{hp04,stdh04}, 
and $(d^2,d)$ ETFs based on the Heisenberg group~\cite{zauner99,rbsc04,grassl04,appleby04}. 
To this list we add a new construction using adjacency matrices of Paley tournaments. 
In section~\ref{sec:construction} we demonstrate the existence of $(2d\minus 1,d)$ and $(2d\minus 1,d\minus 1)$ equiangular tight frames 
whenever $d$ is a power of two or $2d-1\equiv 3\!\!\mod 4$ and is a power of a prime. In section~\ref{sec:connections} we find an 
explicit description of the $(2d\minus 1,d)$ ETFs for odd prime powers by connecting
them to the construction of $(2d,d)$ ETFs using quadratic residues due to Zauner~\cite{zauner99}.


\section{Gram Matrix Construction}
\label{sec:construction}

\noindent To construct an equiangular tight frame, it is sufficient to construct an appropriate Gram matrix for the set. 
\begin{lemma}
\label{lem:gram}
Let $G$ be an $n\times n$ matrix with the following properties:\\
$\phantom{bliiii}$\emph{(i)} $\,|G_{jk}|^2=(n-d)/d(n-1)$, and\\
$\phantom{blah}$\emph{(ii)} $\,G$ has only one nonzero eigenvalue $n/d$ of degeneracy $d$.\\
Then an $(n,d)$ ETF may be explicitly constructed.
\end{lemma}
{\bf Proof}\hspace{.075in}
The frame vectors are given by the columns of the $d\times n$ rank $d$ matrix $T$ satisfying $G=T^* T$~\cite{hj85}. 
More concretely, the frame vectors can be found by factorizing $G$ into $V\Lambda V^*$, where $\Lambda$ is the diagonal matrix of eigenvalues and $V$ is the
unitary matrix whose columns are the eigenvectors. Arranging $\Lambda$ such that the nonzero eigenvalues of $G$ are located in its first $d$ entries,
the $\phi_k$ are given by $\phi_k=\sqrt{n/d}\,\{V_{kl}\}_{l=1}^d$~\cite{sh03}.\hspace{\stretch{1}}$\square$



We shall use adjacency matrices of Paley tournaments in the construction of the Gram matrix. 
Recall that the $n\times n$ adjacency matrix $A$ of the Paley tournament is defined as $A_{jk}=\chi(j\minus k)$ for any prime power 
$n\equiv 3\!\!\mod\!\!$~4. 
Here $\chi(a)$ is the quadratic character of $\mathbb{F}_n$, which is $1$ or $-1$ depending on whether $a$ is or is not a nonzero square and
zero if $a=0$. 
The matrix $A$ has zero diagonal and off-diagonal elements $\pm 1$. It is easily seen by direct calculation that $A^2=J-nI$,
where $J$ is the matrix of all $1$s, using the fact that $\sum_{a\in\mathbb{F}_n}\chi(a)\chi(a-b)=-1$ for all $b\neq 0$~\cite{bjl86}. 
Beyond this case, we also have the following lemma.
\begin{lemma}
\label{lem:conf}
An $n\times n$ zero-diagonal matrix $A$ with off-diagonal elements $\pm 1$ such that $A^2=J-nI$ exists for all $n=2^k-1$, $k\in\mathbb{N}$.
\end{lemma}
\noindent {\bf Proof}\hspace{.075in} Consider the recursive construction of an $(n+1)\times (n+1)$ conference matrix $C$ satisfying $CC^T=nI$:
\begin{eqnarray}
C_2 &=& \left(\begin{array}{cc}0&-1\\1&0\end{array}\right),\\
C_{2m}&=&\left(\begin{array}{cc}C_m&C_m-I_m\\C_m+I_m&-C_m\end{array}\right).
\label{eq:recurse}
\end{eqnarray}
Let $e$ be the column matrix of $n$ $1$s. By induction, $C_{2^k}$ has zero diagonal and the following block form:
\begin{equation}
C_{2^k}=\left(\begin{array}{cc}0 & -e^T\\e&A\end{array}\right).
\end{equation}
By the conference matrix condition, $A$ is the sought-after matrix.\hspace{\stretch{1}}$\square$

Now we can prove our first main result.

\begin{theorem}
Let $n=2d\!-\!1$ for either $n=2^k\!-\!1$, $k\in\mathbb{N}$ or $n=p^k$ for prime $p$ such that $p^k\equiv 3\!\!\mod 4$. 
Then there exist $(n,d)$ and $(n,d\!-\!1)$ equiangular tight frames which can be explicitly constructed. 
\end{theorem}

\noindent {\bf Proof}\hspace{.075in} First we construct the $(n,d)$ ETF. 
The Gram matrix $G$ may be written
\begin{equation}
G_{(n,d)}=\frac{1}{2d}\left(J+nI+i\sqrt{n}A\right).
\end{equation}
Note that because $G$ has complex entries, the resulting frame vectors will
generally be elements of $\mathbb{C}^d$. 
Condition (i) of Lemma~\ref{lem:gram} is clear by inspection. To find the eigenvalues of $G$, note that since $A^2=J-nI$, 
$e$ is an eigenvector of $A$, $J$, and $I$. As $J$ has no other nonzero eigenvalues, 
$A$ and $J+nI$ commute. The eigenvalues are then easily found to satisfy condition (ii).


For the $(n,d\!-\!1)$ equiangular tight frame, take the Gram matrix to be
\begin{equation}
G_{(n,d\!-\!1)}=\frac{1}{2(d\!-\!1)}\left(nI-J+i\sqrt{n}A\right).
\end{equation}
Existence follows from a similar argument to that above.\hfill$\square$


\section{Connections to $(2d,d)$ ETFs}
\label{sec:connections}
Several examples from the previous section were found numerically by starting with $(2d,d)$ ETFs, removing one element,
and forming the canonical tight frame from the resulting set. Recall that for a generic frame consisting of elements $f_k$, the
canonical tight frame formed from $\{f_k\}$ is given by $\{\phi_k=S^{-1/2}f_k\}$ where $S$ is the \emph{frame operator}
defined by $Sf_k=\sum_j\langle f_j,f_k\rangle f_j$~\cite{casazza00}. 

The next two theorems formalize this relationship between the $(2d,d)$ and $(2d\minus 1,d)$ equiangular tight frames for $2d-1\equiv 3\!\mod4$. 
We will make repeated use of the following lemma on sums of quadratic residues, adapted from Chapter 5 of~\cite{ln86}.
\begin{lemma}
\label{lem:gauss}
Suppose $q=p^m$ for prime $p$. For a nontrivial additive character $\psi$ of $\mathbb{F}_q$, define 
\begin{equation}\sigma(a)=\sum_{b\in\mathbb{F}_q}\psi(ab^2).\end{equation}
Then we have for all $a\neq 0$,
\begin{eqnarray} ({\rm i}) &\hspace{.1in}& |\sigma(a)|=\sqrt{q}\;\;\quad\qquad\qquad {\rm and}\\
({\rm ii}) &\hspace{.1in}& \sigma(a)=i\chi(a)(-1)^m\sqrt{q},\quad {\rm for}\,\,q\equiv 3\!\!\mod 4.
\end{eqnarray}
\end{lemma}
\noindent{\bf Proof}\hspace{.075in} Consider the Gaussian sum associated to any additive character $\psi$ and 
multiplicative character $\xi$ of $\mathbb{F}_q$, defined by
\begin{equation}
\Sigma(\psi,\xi)=\sum_{a\in\mathbb{F}_q}\psi(a)\xi(a).
\end{equation}
Call $\psi_1$ the canonical additive character, i.e., $\psi_1(a)=e^{2\pi i{\rm Tr}(a)/p}$ where
Tr is the absolute trace function from $\mathbb{F}_q$ to $\mathbb{F}_p$.
Every additive character is given by some $\psi_c$, where $\psi_c(a)=\psi_1(ca)$. 
Using the quadratic character $\chi$, we have
$\Sigma(\psi_{ba},\chi)=\sum_{c\in\mathbb{F}_q}\psi_{ba}(c^2)=\sum_{c\in\mathbb{F}_q}\psi_b(ac^2)$. 
Since $|\Sigma(\psi,\xi)|=\sqrt{q}$, we have established (i).

For (ii), note that the Gauss sum obeys $\Sigma(\psi_c,\xi)=\overline{\xi(c)}\Sigma(\psi_1,\xi)$, and 
$\Sigma(\psi_1,\chi)=i(-1)^m\sqrt{q}$ for $q\equiv 3\!\!\mod 4$.\hspace{\stretch{1}}$\square$

Now we recapitulate Zauner's construction of $(2d,d)$ ETFs.
\begin{theorem}
A $(2d,d)$ equiangular tight frame exists whenever $q=2d-1$ is any power of an odd prime.
\end{theorem}
\noindent{\bf Proof}\hspace{.075in}
Call the $(q-1)/2$ nonzero 
quadratic residues $b_1,b_2,\dots b_{(q-1)/2}$ and the elements of the field itself $a_1,a_2,\dots a_q$. The columns
of the following $\frac{q+1}{2}\times (q+1)$ matrix $T$ form an ETF:
\begin{equation}
T=\frac{1}{\sqrt{q}}\left(
\begin{array}{ccccc}
\sqrt{q} & 1& 1 &\cdots &1\\
0 & \sqrt{2}\psi(b_1a_1) & \sqrt{2}\psi(b_1a_2) & \cdots & \sqrt{2}\psi(b_1a_q)\\
0 & \sqrt{2}\psi(b_2a_1) & \sqrt{2}\psi(b_2a_2) & \cdots & \sqrt{2}\psi(b_2a_q)\\
\vdots & \vdots & \vdots &\ddots &\vdots\\
0 & \sqrt{2}\psi(b_{\frac{q-1}{2}}a_1) & \sqrt{2}\psi(b_{\frac{q-1}{2}}a_2) & \cdots & \sqrt{2}\psi(b_{\frac{q-1}{2}}a_q)
\end{array}
\right)
\end{equation}
Each column vector $\phi_k$ is manifestly normalized, and the first clearly has overlap $1/\sqrt{q}=1/\sqrt{2d-1}$ with 
any of the others. To establish the remaining pairwise overlaps, it suffices to use Lemma~\ref{lem:gauss}(i).\hspace{\stretch{1}}$\square$


Now follows the connection between these equiangular tight frames and those based on Paley tournaments.
\begin{theorem}
For all $q=p^m=2d\minus 1=3\!\!\mod 4$ a power of a prime $p$, there exists a $(2d-1,d)$ ETF with Gram matrix 
$G=(J+qI+i\sqrt{q}A)/2d$ which is given by 
the canonical tight frame associated with the set of vectors formed by removing the first element from 
Zauner's construction. 
\end{theorem}

By removing the first column from the matrix $T$ and forming the frame operator $S=TT^*$, we find that 
\begin{eqnarray*}
S_{00} &=& 1,\\
S_{j0}&=&S_{0k}^* = \frac{\sqrt{2}}{q}\sum_m\psi(b_la_m)=0,\\
S_{jk}&=&\frac{2}{q}\sum_m \psi((b_j-b_k)a_m)=2\delta_{jk},
\end{eqnarray*}
so $S={\rm diag}\{1,2,\dots,2\}$.
The new matrix $\widetilde{T}$, whose columns are the normalized vectors $\sqrt{q/d}\,S^{-1/2}\phi_k$, is written as
\begin{equation}
\widetilde{T}=\sqrt{\frac{1}{d}}\left(
\begin{array}{cccc}
 1& 1 &\cdots &1\\
 \psi(b_1a_1) & \psi(b_1a_2) & \cdots & \psi(b_1a_q)\\
 \psi(b_2a_1) & \psi(b_2a_2) & \cdots & \psi(b_2a_q)\\
 \vdots & \vdots &\ddots &\vdots\\
 \psi(b_{\frac{q-1}{2}}a_1) & \psi(b_{\frac{q-1}{2}}a_2) & \cdots & \psi(b_{\frac{q-1}{2}}a_q)
\end{array}
\right).
\end{equation}

The Gram matrix is easily found using Lemma~\ref{lem:gauss}(ii).\hspace{\stretch{1}}$\square$



We expect that a similar relationship holds between the $(2d,d)$ and $(2d-1,d)$ ETFs constructed via
the Conference matrices above, though at present we have only numerical evidence to support this claim. 

The author thanks K.~Manne and C.~M.~Caves for helpful comments on this manuscript. This work was
supported in part by the Office of Naval Research Grant No. N00014-03-1-0426.


\begin{thebibliography}{99}

\bibitem{bf03}J.J.~Benedetto, M.~Fickus, Finite 
normalized tight frames, Adv. Comp. Math. 18 (2003) 357--385.

\bibitem{sh03}T.~Strohmer, R.~W.~Heath, Jr., Grassmannian 
frames with applications to coding and communication, 
Appl. Comp. Harm. Anal. 14 (2003) 257--275.


\bibitem{Koenig03}
H.~K\"onig, Cubature formulas on spheres, in: W.~Haussmann, K.~Jetter, M.~Reimer (Eds.), Advances in Multivariate Approximation: Proceedings of the 3rd 
International Conference on Multivariate Approximation Theory, vol.~107 of Math. Res., Wiley-VCH, Berlin, 1999, pp.~201--211.

\bibitem{zauner99}
G.~Zauner, Quantum Designs --- Foundations of a Non-Commutative Theory of Designs (German),
PhD thesis, University of Vienna, 1999. 


\bibitem{ls73}{ P.~W.~H.~Lemmens, J.~J.~Seidel}, Equi-isoclinic subspaces of Euclidean space, Konink. Ned. Akad. Wet. A 76 (1973)  98--107.


\bibitem{chs96} { J.~H.~Conway, R.~H.~Hardin, N.~J.~A.~Sloane}, 
Packing lines, planes, etc.: Packings in Grassmannian spaces, 
J. Exp. Math. 5 (1996) 139--159.

\bibitem{hp04}R.~B.~Holmes, V.~I.~Paulsen, Optimal frames for erasures, Lin. Alg. Appl. 377 (2004) 31--51.

\bibitem{stdh04}{ M.~A.~Sustik, J.~A.~Tropp, I.~S.~Dhillon, R.~W.~Heath, Jr.}, On the existence of equiangular tight frames,
submitted to Lin. Alg. Appl. June 2004. 

\bibitem{rbsc04}
J.~M.~Renes, R.~Blume-Kohout, A.~J.~Scott, C.~M.~Caves, 
Symmetric informationally complete quantum measurements,
J.~Math.~Phys. 45 (2004) 2171--2180.
b
\bibitem{grassl04}
M.~Grassl, On SIC-POVMs and MUBs in dimension 6, quant-ph/0406175.

\bibitem{appleby04}
D.~M.~Appleby, SIC-POVMs and the Extended Clifford Group, quant-ph/0412001.

\bibitem{hj85}R.~A.~Horn, C.~R.~Johnson, Matrix Analysis, Cambridge University Press, 1985.

\bibitem{bjl86} T.~Beth, D.~Jungnickel, H.~Lenz, Design Theory, Cambridge University Press, 1986.

\bibitem{casazza00}P.~G.~Casazza, The art of frame theory, Taiwan. J. Math. 4 (2000) 129--201.

\bibitem{ln86}{ R.~Lidl, H.~Niederreiter}, Introduction to Finite Fields and their Applications, Cambridge University Press, 1986. 


\end{thebibliography}
\end{document}